   \documentclass[12pt,onecolumn,twoside]{IEEEtran}

\IEEEoverridecommandlockouts              % This command is only needed if 
\usepackage{graphicx}

\usepackage{url}
\usepackage{amsmath,amssymb}
\usepackage{amsthm,thmtools}
\usepackage{amsfonts} 

\usepackage{enumerate}
\usepackage{color}                  %******** for colored remarks********

\usepackage{verbatim}

\usepackage{amssymb}
\usepackage{amsbsy}
\usepackage{amsmath}
\usepackage{setspace}
\usepackage{url}

\usepackage{float}
 \usepackage{dsfont}

 \newcommand{\diag}{\operatorname{diag}}
 
 \newcommand{\tr}{\operatorname{tr}}

\newcommand{\st}{\, | \,}

\newcommand{\spanop}{\operatorname{span}}

\newcommand*\diff{\mathop{}\!\mathrm{d}}

\declaretheorem[name={Example},qed={\lower-0.3ex\hbox{$\square$}} ] {Example}

\declaretheorem[name={Definition}  ] {Definition}
\declaretheorem[name={Theorem}  ] {Theorem}
\declaretheorem[name={Lemma}  ] {Lemma}
\declaretheorem[name={Remark}  ] {Remark}

\declaretheorem[name={Proposition}  ] {Proposition}
\declaretheorem[name={Fact}  ] {Fact}

\newcommand {\R}{\mathbb R}

\newcommand{\be}{\begin{equation}}
\newcommand{\ee}{\end{equation}}

\usepackage{lineno}
\begin{document}

%
%\onecolumn  % add this line after \begin{document} but before \titlepage

%\doublespace
\title{From   Partial and Horizontal Contraction to $k$-Contraction}
%\title{Relations between partial contraction, horizontal contraction, and $k$-contraction
%\thanks{This research is
%		supported in part 
%		by Knut and Alice Wallenberg Foundation
% by the Swedish Research Council, research grants  from  the Israel Science Foundation  and the 
%		US-Israel Binational Science Foundation. \emph{(Corresponding author: Michael Margaliot.) }
%}}

\author{Chengshuai Wu$^*$ and Dimos V. Dimarogonas
	\thanks{
	\IEEEcompsocthanksitem
	This work is supported by Knut and Alice Wallenberg Foundation and the Swedish Research Council. 
	    \IEEEcompsocthanksitem
		C. Wu and D. V. Dimarogonas are with the School of Electrical Engineering and Computer Science, KTH Royal Institute of Technology and Digital Futures, Stockholm, SE-100 44, Sweden (\texttt{chengshuai.wu@gmail.com}, \texttt{dimos@kth.se}).
		\IEEEcompsocthanksitem
		$^*$ Corresponding author.
		}
}

\maketitle
%\doublespace 
%%%%%%%%%%%%%%%%%%%%%%%%%%%%%%%%%%%%%%%%%%%%%%%%%%%%%%%%%
 %\begin{multicols}{2} 
%\begin{center} 
					%			\today \\
%\end{center}

\begin{abstract}
%%%%%%%%%%%%%%%%%
A geometric generalization of contraction theory called~$k$-contraction was recently developed using $k$-compound matrices. 
In this note, we focus on the relations between $k$-contraction and two other generalized contraction frameworks:  partial contraction (also known as virtual contraction) and horizontal contraction. We show that in general these three notions of contraction are different. We here provide new sufficient conditions guaranteeing that partial contraction implies horizontal contraction, and that horizontal contraction implies $k$-contraction. We use the Andronov-Hopf oscillator to demonstrate some of the theoretical results.
%%%%%%%%%%%%%%%%%
\end{abstract}
%%%%%%%%%%%%%%%%%%%%%%%%%%%%%%%%%%%%%%%%%%%%%%%%%%%%%%%%%%%%%%%%%%%%%%%%%%%%


\section{Introduction} \label{sec:intro}  
Contraction theory is a powerful tool for analyzing the asymptotic behavior of 
nonlinear time-varying dynamical systems \cite{LOHMILLER1998683, sontag_cotraction_tutorial, forni2014}. A contractive  system
behaves in many respects like
a uniformly asymptotically stable linear system:  initial conditions are ``forgetten'' and 
any two trajectories approach each other at an exponential rate.
%%%%

There exist easy to verify sufficient conditions for contraction that are based on  matrix measures \cite{sontag_cotraction_tutorial} and contraction analysis has found numerous applications such as control synthesis for regulation~\cite{pavlov2008incremental}
and tracking~\cite{wu2019robust},  
observer design \cite{doi:10.1002/aic.690460317,sanfelice2011convergence, aghannan2003intrinsic}, optimization \cite{wensing2020beyond}, synchronization of multi-agents systems \cite{slotine2005study,russo2009solving}, robotics \cite{manchester2018unifying}, learning algorithm~\cite{wensing2020beyond},
and systems biology~\cite{RFM_entrain,entrain2011}. 

Any two solutions of a contractive system converge to each other, which implies a unique exponentially asymptotically stable equilibrium or trajectory. This rules out the existence of multiple (stable or unstable) equilibriums, limit cycles, and other oscillatory behaviors. This motivates researchers to introduce generalizations of contraction theory which allow analyzing non-trivial attractors, for example, \emph{partial contraction} \cite{wang2005partial,belabbas2010factorizations,slotine2005study}, \emph{horizontal contraction} \cite{forni2014}, and $k$-contraction \cite{kordercont}. 

Roughly speaking, partial contraction is related to the contractive behavior of an auxiliary system associated with the studied one, and horizontal contraction studies contractive properties along some particular ``directions''. Despite using different mathematical formulations, these two generalized notions of contraction are both effective for analyzing stable limit cycles or synchronization of networked systems, i.e., convergence to certain subspaces. 
%%%%%%%
Furthermore, Ref.~\cite{wang2016immersion} showed that both notions (where partial contraction is referred to as virtual contraction) can be utilized to solve a particular control problem arising in the \emph{immersion $\&$ invariance} (I$\&$I) stabilization procedure \cite{astolfi2003immersion}.
%%%%%%%%%%%

Ref.~\cite{kordercont}
introduced a generalization of contraction theory called $k$-contraction  based on the seminal work of Muldowney~\cite{muldo1990}. A dynamical system is called $k$-contractive if the dynamics contracts $k$-volumes at an exponential rate. For~$k=1$, this recovers the standard contraction theory as $1$-volume is just length. 
However, $k$-contraction with~$k>1$ can be applied to analyze systems that are not contractive (i.e., not $1$-contractive) such as multi-stable systems that are prevalent in mathematical models of real-world systems.  
In particular, it was  shown in \cite{kordercont,wu2022generalization} that $k$-contraction can also be applied to chaotic systems, which typically cannot be analyzed using
partial/horizontal contraction.

\begin{figure}
	\begin{center}
		\centering
		\includegraphics[scale=0.35]{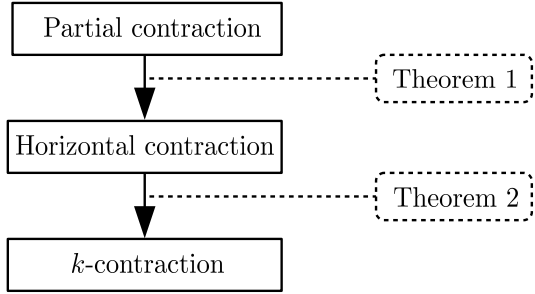} 
		\caption{Illustration of the main results.}
		\label{block} 
	\end{center}
\end{figure}

\begin{comment}
Since the Hausdorff dimension of the attractors of $k$-contractive systems is less than $k$, more general attractors are possible.
\end{comment}
However, the three notions are closely related in certain cases and it is our intention to bring more insights into the distinctions and the relations. 
 %%%%
Unlike \cite{kordercont}, this current work considers dynamical systems whose solutions evolve on a forward invariant and connected, but not necessarily convex manifold.
%%%
In our main results, we provide conditions describing when partial contraction implies horizontal contraction and when horizontal contraction implies $k$-contraction.  
By combining these results together, sufficient conditions for partial contraction implying $k$-contraction are also obtained (see Fig.~\ref{block}).
%%%
Furthermore, some examples are given to show that the reverse implications do not hold in general. These results are useful since $k$-contraction and, in particular, $2$-contraction implies strong results on the attractors and asymptotic behavior of nonlinear time-invariant systems~\cite{li1995,kordercont}. Therefore, the same conclusions can be drawn for partially or horizontally contractive systems if the aforementioned sufficient conditions hold.

The remainder of this note is organized as follows. The next section briefly reviews the definitions of $k$-contraction, horizontal contraction, and partial contraction. Section \ref{sec:main1} and \ref{sec:main2} detail the main results.
The Andronov-Hopf oscillator is revisited in Section \ref{sec:aho} to validate the proposed results. Section \ref{sec:con} gives the conclusions.

\textbf{Notation}.
%%%%%%%%%%%%%%%%%%%%%%%%
%\section{Preliminaries} \label{sec:pre}
%%%%%%%%%%%%%%%%%%%%%%%%%%%%% 
Here we briefly describe the basic notations and some mathematical tools including compound matrices and wedge products are required to define $k$-contraction, see~\cite{comp_long_tutorial,muldo1990,fiedler_book} for more details and proofs. 

For two integers~$i,j$, with~$0<i\leq j$, we denote~$[i,j]:=\{i,i+1,\dots,j\}$.
%%%%%%%%%%%%%%
For an $n$-dimensional manifold $\Omega \subseteq \R^n$, we denote the \emph{tangent space} of $\Omega$ at $x \in \Omega$ by $T_x \Omega$, and the \emph{tangent bundle} of $\Omega$ by $T\Omega := \cup_{x \in \Omega} (\{x\} \times T_x \Omega)$. 

%%%%%%%%%%%%%%%%%%%%%%%%%%%%%%%%%%%%%%%%%%%%
%\subsection{Compound matrices}\label{sec:compound}
%%%%%%%%%%%%%%%%%%%%

\emph{Compound matrices}.
\begin{comment}
 For $k \in [1,n]$, let~$Q_{k,n}$ denote the set 
of increasing sequences of~$k$ numbers from~$[1,n]$
ordered lexicographically. With a slight abuse of notation, we will sometimes treat such ordered sequences as sets. The cardinality of $Q_{k,n}$ is $r:=\binom{n}{k}$. The $i$th element of $Q_{k,n}$ is denoted $\kappa_i$. For example,
$
Q_{2,3} =\{\kappa_1, \kappa_2, \kappa_3\}, \text{ with } ~\kappa_1= \{1,2\},~\kappa_2 = \{1,3\}, ~\kappa_3 = \{2,3\}.
$

Given~$A\in\R^{n\times m}$ and~$k\in[1,\min\{n,m\}]$, 
a minor of order~$k$ of~$A$ is the determinant of some~$k \times k$
submatrix of~$A$.  
%%%%%%%%%
%%%%%%%%%
Consider the
$\binom{n}{k} \binom{m}{k}  $
 minors of  order~$k$ of~$A$. 
Each such minor is defined by a set of row indices~$\kappa_i \in Q_{k,n}$ and column indices~$\kappa_j \in 
Q_{k,m}$. This minor 
is denoted by~$A(\kappa_i|\kappa_j)$.
For example, for~$A=\begin{bmatrix} 4&5&8   \\ -1 &4 &9  \\0&3 &7
\end{bmatrix}$, one of its minors of order $2$ is
\[
A(\{ 1,3\} |\{1,2\})=\det  \begin{bmatrix} 4&5\\0&3
\end{bmatrix}   =12.
\]
\end{comment}
%%%%%%%%%%%
Given~$A\in\R^{n\times m}$ and~$k\in[1,\min\{n,m\}]$, the~$k$th \emph{multiplicative compound matrix} 
of~$A$, denoted~$A^{(k)}$, is the~$\binom{n}{k}\times  \binom{m}{k}$ matrix
that includes all the minors of order~$k$ of $A$ ordered lexicographically. 
\begin{comment}
That is, the $ij$th entry of $A^{(k)}$ is $A(\kappa_i \vert \kappa_j)$.
\end{comment}
%%%%%%%%%%%%%%%%%%%%%%%%%%%%%%%%%%%%%%%%%%%%%%
In particular, $A^{(1)}=A$, and for $A \in \R^{n \times n}$, $A^{(n)}=\det(A)$. 
\begin{comment}
In particular, for the identity matrix $I \in \R^{n \times n}$, $I^{(k)}$ is the~$r\times r$ identity matrix, with~$r:=\binom{n}{k}$. 
\end{comment}
The Cauchy-Binet formula  (see, e.g.,~\cite[Ch.~1]{total_book}) asserts   
that for any~$A \in \R^{n \times p}$, $B \in \R^{p \times m}$, and $k = [1,\min\{n,m,p\}]$,
\be \label{eq:prodab}
%%%
(AB)^{(k)}=A^{(k)} B^{(k)},
\ee
which justifies the term multiplicative compound. 

\begin{comment}
When~$n=p=m=k$ this becomes the familiar formula~$\det(AB)=\det(A)\det(B)$.
If~$A$ is square and non-singular then~\eqref{eq:prodab} implies that~$ I^{(k)}=(AA^{-1})^{(k)}=A^{(k)} (A^{-1})^{(k)}$,
	so~$(A^{(k)})^{-1}=(A^{-1})^{(k)} $. 
 %%%%%%%%%%%%%%%%%%%%%%%%%%%%%%%%%
\end{comment}

The $k$th \emph{additive compound matrix} of a square matrix~$A \in \R^{n \times n}$
is defined  by
\[
	A^{[k]}:= \frac{d}{d \varepsilon}  (I_n +\varepsilon A)^{(k)} |_{\varepsilon=0},
\]
where $I_n$ denotes the $n \times n$ identity matrix. 
This implies that
\be\label{eq:poyrt}
(I_n +\varepsilon A)^{(k)}= I_r +\varepsilon   A^{[k]} +o(\varepsilon ) ,
\ee 
%%%%%%
%%%
with~$ r:= \binom{n}{k}$. 
The matrix~$A^{[k]}$ can be given explicitly in terms of the entries~$a_{ij}$ of~$A$ as shown in \cite{schwarz1970, fiedler_book}.  
It follows from~\eqref{eq:poyrt} and the properties of the multiplicative compound that for any $A, B \in \R^{n \times n}$
\be \label{eq:addi}
(A+B)^{[k]}= A^{[k]}+B^{[k]},
\ee
which justifies the term additive compound.

 \emph{Wedge products}.
%%%%%%%%%%%%%%%%%%%%%%%%%%%%%%%%%%%%%%%%%%%%%
%\subsection{Wedge products} \label{sec:wedge}
%%%%%%%%%%%%%%%%%%%%%%%%%%%%%%%%%%%%%%%%%%%%%%
The multiplicative compound matrix has an important  geometric interpretation in terms of the $k$-volume 
of a  $k$-parallelotope~\cite{comp_long_tutorial}.  Pick~$k \in [1,n]$, and~$k$ vectors~$a^i \in \R^n$, $i =  1, \dots, k$. The wedge product of these vectors, denoted
$	a^1\wedge\dots\wedge a^k$, can be represented using the multiplicative compound as 
\be\label{eq:wedgdef}
	a^1\wedge\dots\wedge a^k=		\begin{bmatrix} a^1 &\dots & a^k \end{bmatrix}^{(k)}.
\ee
This provides a representation of~$a^1\wedge\dots\wedge a^k$ as an~$r$-dimensional column vector, where~$r:=\binom{n}{k}$. We will  use the short-hand notation~$\wedge_{i=1}^k a^i := a^1\wedge\dots\wedge a^k$ throughout the paper.
%%%%%%%%%%%
The $k$-parallelotope generated by $a^1, \dots, a^k$ (and the zero vertex) is~$P(a^1,\dots,a^k):=\{ \sum_{i=1}^k 
c_i a^i \st 
 c_i\in[0,1]\} $.
The $k$-volume of  $P(a^1, \dots, a^k)$
is $|\wedge_{i=1}^k a^i |_2$, where $|\cdot|_2$ is the~$L_2$ norm. 
In the particular case~$k=n$ this reduces to the well-known formula
\[
\text{volume}(P(a^1,\dots,a^n))=|\det(a^1,\dots, a^n)|.
\]
\section{Three Generalized Notions of Contraction}
\label{sec:three}
%% \section{Preliminaries}
In this section, we briefly review the definitions of $k$-contraction,  horizontal contraction, and partial contraction. 

Consider the nonlinear time-varying (NTV) system
%%%%%%%%%%%%%%%%%%%%%%
\be\label{eq:nonlinsys}
\dot x = f(t,x),
\ee
where $f :\R_+\times \R^n \to \R^n$ is continuously differentiable w.r.t. its second argument, and let
\be \label{eq:jaco}
J(t,x):=\frac{\partial}{\partial x}f(t,x)
\ee
denote the Jacobian of the vector field w.r.t.~$x$.
%%%%%%%%%%
We denote by~$x(t, t_0, a)$ the solution to \eqref{eq:nonlinsys} at time~$t$ emanating from the initial condition~$a \in \R^n$ at time $t_0$, that is, $x(t_0, t_0, a) = a$. 
We assume that the solutions of \eqref{eq:nonlinsys} evolve on a closed and connected $n$-dimensional manifold~$\Omega$, and that for any initial condition~$a \in \Omega$, a unique solution~$x(t, t_0, a)$ exists and satisfies~$x(t,t_0,a)\in\Omega$ for all~$t\geq t_0$. For the sake of simplicity, we assume from here on
that the initial time is~$t_0 = 0$, and write $x(t, a) := x(t, 0, a)$.
%When clear from the context, we sometimes use $x$ or $x(t)$ to denote $x(t, a)$.

Consider the matrix 
\be
\Phi(t,a) := \frac{\partial x(t,a)}{\partial a}.
\ee
Note that $\Phi(0,a) = I_n$. A straightforward computation yields
\be \label{eq:ptrma}
%%%
\frac{d}{dt} \Phi(t,a) =J(t, x(t,a)) \Phi(t,a) . 
%%%%
\ee
This is  the \emph{variational system} associated with \eqref{eq:nonlinsys} along $x(t,a)$. 
%%%%%%%%%%%
Let $\delta a \in T_a \Omega$ denote an infinitesimal   variation to the initial condition~$a \in \Omega$.
Then~$\delta x(t,a) := \Phi(t,a)\delta a$ is the infinitesimal displacement w.r.t. the solution $x(t,a)$ induced by the initial condition~$a + \delta a$.  Eq.~\eqref{eq:ptrma} implies that 
%%%%%%%%%%%%%%%%%%%
\be \label{eq:vardx}
\delta \dot{x}(t,a) := \frac{\diff}{\diff t} \delta x(t,a) = J(t, x(t,a)) \delta x(t,a).
\ee

%%%%%%%
\subsection{$k$-contraction} \label{sec:koc}
For $k \in [1,n]$, consider~$\Phi^{(k)}(t,a):=(\Phi (t,a))^{(k)}$, i.e., the $k$th multiplicative compound matrix of $\Phi (t,a)$. Fix~$\varepsilon>0$. Eq.~\eqref{eq:ptrma} and the Cauchy–Binet formula give 
\begin{align*}
      \Phi^{(k)} (t+\varepsilon,a)&=  (\Phi(t,a)+\varepsilon J(t, x(t,a))\Phi(t,a)   ) ^{(k)}+o(\varepsilon)\\
			&=(I_n+\varepsilon J(t, x(t,a))   ) ^{(k)} \Phi  ^{(k)}(t,a)+o(\varepsilon).
\end{align*}
Combining this with~\eqref{eq:poyrt} and the fact that~$\Phi(0,a)=I_n$ yields  a differential equation for $\Phi^{(k)}(t,a)$:
\be\label{eq:povt}
\frac{d}{dt} \Phi^{(k)}(t,a)=J^{[k]}(t, x(t,a)) \Phi^{(k)}(t,a), ~ \Phi^{(k)}(0,a)=I_r, 
\ee
where~$ J^{[k]}(t, x(t,a)):=(J(t, x(t,a))) ^{[k]}$, and $r: = \binom{n}{k}$. 
In other words, all the minors of order~$k$ of~$\Phi(t,a)$, stacked in the matrix~$\Phi^{(k)}(t,a)$, satisfy   
a linear dynamics with the matrix~$J^{[k]}(t, x(t,a))$. 

Pick $k$ initial conditions~$\delta a^1, \dots, \delta a^k$ for \eqref{eq:vardx}. Define 
\be \label{eq:defy}
\delta x^i(t, a) := \Phi(t,a) \delta a^i, \quad y(t,a): = \wedge_{i=1}^k \delta x^i(t,a) .
\ee
Note that $y(t,a) = 0$ iff $\delta x^1(t,a), \dots, \delta x^k(t,a)$ are linearly dependent. By \eqref{eq:prodab} and \eqref{eq:wedgdef}, 
\be \label{eq:wedgey}  \begin{aligned}
y(t,a) = \wedge_{i =1}^k \Phi(t,a) \delta a^i 
  = \Phi^{(k)}(t,a) y(0,a),
\end{aligned}\ee
%%%%%%%%%
\begin{comment}
=& \begin{bmatrix}
\Phi(t,a) \delta a^1 & \dots &  \Phi(t,a) \delta a^k  
\end{bmatrix}^{(k)} \\
=& \Phi^{(k)}(t,a) \begin{bmatrix}
\delta a^1 & \dots &   \delta a^k  
\end{bmatrix}^{(k)} \\
\end{comment}
%%%%%%%%%%%%
and \eqref{eq:povt}   yields
\be \label{eq:dyta}
\dot{y}(t,a) = J^{[k]}(t, x(t,a)) y(t,a).
\ee
This is the \emph{$k$th compound equation} of \eqref{eq:vardx} along $x(t,a)$ (see e.g.~\cite{LI2000295}).
%%%%%%%%%%%%%
This leads to the  following  definition.
%%%%%%%%%%%%%%%%%
\begin{Definition}[$k$-contraction] \label{def:kord}
Fix   $k \in [1,n]$ and let $r:=\binom{n}{k}$. The NTV system \eqref{eq:nonlinsys} is called \emph{$k$-contractive} if the linear time-varying (LTV) system
	\begin{equation}
	\dot{y}(t) = J^{[k]}(t, x(t,a)) y(t)
	\end{equation}
	is uniformly exponentially stable for any $a \in \Omega$, that is, there exist~$c \geq 1$, $\eta > 0$, and a vector norm $\vert \cdot \vert$ such that
	\be \label{eq:def}
	\vert y(t) \vert \leq c e^{-\eta t} \vert y(0) \vert, \text{ for all } ~t \geq 0.
	\ee
\end{Definition}

Note that the above definition is slightly different from \cite[Def. 2]{kordercont} where the value of $c$ is fixed as one. 
%%%%%%%%%%%%%%%%%%%
By the geometric interpretation of wedge products, Eq.~\eqref{eq:def} implies that  the $k$-volume of the $k$-parallelotope generated by the vertices~$\delta x^1(t,a), \dots, \delta x^k(t,a)$ (and the zero vertex)  decays to zero at an exponential rate.
 For $k =1$, $k$-contraction 
 reduces to   standard contraction.

\subsection{Horizontal contraction} \label{sec:hori}
 For every $x \in \Omega$, suppose that $T_x \Omega$ can be subdivided into a horizontal distribution $\mathcal{H}_x$ and a vertical distribution $\mathcal{Q}_x$ which are orthogonally complementary to each other. That is, there exist  $\ell \in [1, n]$ and~$C^1$ mappings $h^i, q^i: \R^n \to  \R^n$ such that
\be \begin{aligned} \label{eq:hqx}
	& \mathcal{H}_x := \spanop\{h^1(x), \dots, h^\ell(x) \},  \\
	& \mathcal{Q}_x := \spanop\{q^1(x), \dots, q^{n-\ell}(x) \} .
\end{aligned}\ee
 Note that if $\ell = n$, then $\mathcal{H}_x = T_x \Omega$. 
Define the matrices
\begin{equation} \label{eq:hqmatrix} \begin{aligned}
& H(x) : = \begin{bmatrix} h^1(x) & \cdots & h^\ell(x) \end{bmatrix} \in \R^{n \times \ell}, \\
& Q(x) : = \begin{bmatrix} q^1(x) & \cdots & q^{n-\ell}(x) \end{bmatrix} \in \R^{n \times (n -\ell)}.
\end{aligned}\end{equation}
Since $\mathcal{H}_x$ and $\mathcal{Q}_x$ are orthogonal to each other, we have 
\be \label{eq:htq0}
H^T(x) Q(x) = 0.
\ee

%%%%%%%%%%%
For every $\delta x \in T_x \Omega$, there exists a set of uniquely defined $\delta x_h \in \R^\ell$ and $\delta x_q \in \R^{n -\ell}$ such that 
\begin{equation} \label{eq:deltahq}
\delta x  = H(x) \delta x_h + Q(x) \delta x_q.
\end{equation}
Note that $H(x) \delta x_h \in \mathcal{H}_x$, and $ Q(x) \delta x_q \in \mathcal{Q}_x$. Combining~\eqref{eq:deltahq} and~\eqref{eq:htq0} gives 
\be \label{eq:hty} \begin{aligned}
& H^T(x) \delta x =  H^T(x) H(x) \delta x_h, \\
& Q^T(x) \delta x =  Q^T(x) Q(x) \delta x_q, 
 \text{ for all } (x, \delta x) \in T \Omega.
\end{aligned}\ee
Without loss of generality, we assume throughout that both $H(x)$ and $Q(x)$ are bounded on $x \in \Omega$. 
%%%%%%%%%%%%%%
Based on the above discussions, horizontal contraction is defined as follows.

\begin{Definition}[Horizontal contraction] \label{def:hori}
The NTV system \eqref{eq:nonlinsys} is called \emph{horizontally contractive} w.r.t. $\mathcal{H}_x$ if there exist $c \geq 1$, $\eta > 0$, and a vector norm $|\cdot|$ such that the solution of \eqref{eq:vardx} for any $a \in \Omega$, i.e., $\delta x(t,a)  = H(x(t,a)) \delta x_h(t,a) + Q(x(t,a)) \delta x_q(t,a)$, satisfies
\be \label{eq:d12345}
\vert \delta x_h(t,a) \vert \leq  c e^{-\eta t}  \vert \delta x_h(0,a) \vert, \quad \text{for all } t\geq 0.
\ee
\end{Definition}

In \cite{forni2014}, horizontal contraction is formalized via a differential Lyapunov framework. Specifically, a sufficient condition is given in terms of a so-called \emph{Horizontal Finsler-Lyapunov function}.  

\begin{Definition}[Horizontal Finsler-Lyapunov function \cite{forni2014}] \label{def:semifins}
	Consider a manifold $\Omega$ and the tangent space $T_x \Omega = \mathcal{H}_x \oplus \mathcal{Q}_x$, where $\oplus$ denotes the direct sum of vector spaces.
	A $C^1$ function $V: T \Omega \to \R_+$ is called a candidate  \emph{horizontal Finsler-Lyapunov function} for \eqref{eq:nonlinsys} if there exist constants $d_1, d_2>0$, $d_3>1$, and a function $F: T \Omega \to \R_+$ such that
	\be\begin{aligned} \label{eq:d123}
		&V(x, \delta x) = V(x, H(x) \delta x_h),\\
		& d_1 (F(x, \delta x))^{d_3} \leq V(x, \delta x) \leq d_2 (F(x, \delta x))^{d_3}, \\
		& \quad \quad \quad \quad \quad \quad \quad \quad \text{ for all } (x, \delta x) \in T \Omega,
	\end{aligned}\ee 
	and $F$ satisfies the following conditions:
	\begin{enumerate}[(i)]
		\item \label{pp:h2} $F(x, \delta x) = F(x, H(x)\delta x_h)$ for every $ (x,\delta x) \in T_x \Omega$;
		\item \label{pp:h1} $F(x,\delta x)$ is $C^1$ for all $x \in \Omega$ and  $ \delta x \in \mathcal{H}_x \setminus \{ 0\}$;
		\item \label{pp:h3} $F(x, \delta x) \geq 0$ for all $(x, \delta x) \in T \Omega$ with equality only when $\delta x \in \mathcal{Q}_x$;
		\item \label{pp:h4}  $F(x,\lambda \delta x) = \vert \lambda \vert F(x, \delta x)$ for all $(x, \delta x) \in T \Omega$ and any $\lambda \in \R$;
		\item \label{pp:h5} $F(x, \delta x^1 + \delta x^2) \leq F(x, \delta x^1) + F(x, \delta x^2) $ for all $(x, \delta x^1),(x,\delta x^2) \in T \Omega$, with equality only when $ \delta x^1_h = \lambda \delta x^2_h$ for some $\lambda \in \R$;
		\item \label{pp:d45} there exist constants $d_4, d_5 >0$ and a vector norm $\vert \cdot \vert$ such that $d_4 \vert \delta x_h \vert \leq F(x,\delta x) \leq d_5 \vert \delta x_h \vert$ for all $(x,\delta x) \in T \Omega$.
	\end{enumerate}  
\end{Definition}

\begin{Proposition}\cite{forni2014} \label{prop:hori}
	Consider the NTV system \eqref{eq:nonlinsys} and the associated variational system \eqref{eq:vardx} with a candidate horizontal Finsler-Lyapunov function $V: T \Omega \to \R_+$. If there exists $\lambda > 0$ such that
	\be \label{eq:dvx}
	\frac{\partial V(x, \delta x)}{\partial x} f(t,x) + \frac{\partial V(x, \delta x)}{\partial \delta x} J(t, x) \delta x \leq -\lambda V(x, \delta x), 
 	\ee
	for all $t \in \R_+$ and all $(x,\delta x) \in T \Omega$, then the system \eqref{eq:nonlinsys} is \emph{horizontally contractive} w.r.t. $\mathcal{H}_x$. 
\end{Proposition}

By taking the derivative of $V(x(t), \delta x(t))$ along the trajectories of \eqref{eq:nonlinsys} and \eqref{eq:vardx}, we have 
\be \label{eq:ddvx}
\dot{V}(x , \delta x ) \leq -\lambda V(x, \delta x), \quad \text{for all } t\geq 0, 
\ee
so $V(x(t), \delta x(t)) \leq
\exp(-\lambda t )  V(x(0), \delta x(0))$. By \eqref{eq:d123} and Property \eqref{pp:d45} in Definition~\ref{def:semifins}, Eq. \eqref{eq:d12345} holds with $c := \left ( {d_2}/{d_1} \right )^{\frac{1}{d_3}}   d_5 / d_4 $ and $\eta  :=   \lambda/ d_3$. 
%%%%%%%%%%

If $\mathcal{Q}_x = \{ 0 \}$ in Definition~\ref{def:hori}, i.e., $\mathcal{H}_x = T_x \Omega$ and $H(x) = I_n$, then \eqref{eq:d12345} is the same as \eqref{eq:def} in Definition~\ref{def:kord} with $k=1$, since $y(t)$ in \eqref{eq:def} is a solution of the variational system \eqref{eq:vardx} in this case. Therefore, Definition~\ref{def:hori} with $\mathcal{H}_x = T_x \Omega$ reduces to the definition for standard contraction. In this case, $V(x, \delta x)$ in Definition~\ref{def:semifins} is called a \emph{Finsler-Lyapunov function} \cite[Def. 2]{forni2014}.

\subsection{Partial contraction}
The following definition is a time-varying version of partial contraction as given in \cite{belabbas2010factorizations}.

\begin{Definition}[Partial contraction] \label{def:parcon}
	Consider the NTV system~\eqref{eq:nonlinsys}. Assume that there exist functions $p(x): \R^n \to \R^\ell$, with~$\ell \in [1,n]$,  and $g: \R_+ \times \R^\ell \times \R^n \to \R^n$ such that 
	\be \label{eq:pxf}
	g(t, p(x), x) =  f(t,x).
	\ee  
	System \eqref{eq:nonlinsys} is called partially contractive w.r.t. $p(x)$ if the system 
	\be \label{eq:xisy}
	\dot{\xi} = f_\xi (t, \xi, x)  := \frac{\partial p(x)}{\partial x} g(t,\xi, x)
	\ee
	is contractive w.r.t. $\xi$ for all $t\in \R_+$, $x \in \Omega$, and $\xi \in \Omega_\xi$. Here, $\Omega_\xi \subseteq \R^{\ell}$ denotes the state-space of \eqref{eq:xisy}. 
\end{Definition}

Here, the functions $p(\cdot)$ and $g(\cdot)$ are related to a \emph{factorization} of $f(t,x)$ \cite{belabbas2010factorizations}. Partial contraction implies that every solution of \eqref{eq:xisy} converges to $p(x(t))$ exponentially since~$\xi(t) =  p(x(t)) $ is a particular solution of~\eqref{eq:xisy}. For the special case that $p(x) = x$, the system~\eqref{eq:xisy} may serve as an observer for~\eqref{eq:nonlinsys} \cite{wang2005partial}.  

\begin{Remark} \label{re:p-m}
Similar to \cite{belabbas2010factorizations}, if we assume that there exists a mapping $m: \R^n \to \R^\ell$ such that 
\be
\frac{\partial p(x)}{\partial x}g(t, m(x), x) = \frac{\partial m(x)}{\partial x} g(t, p(x), x), 
\ee
for all $t \in \R_+$ and $x \in \Omega$, then by \eqref{eq:pxf} and \eqref{eq:xisy}, we have 
\be \label{eq:p-m}
\dot{m}(x) =  f_\xi(t, m(x),x).
\ee
That is, $\xi(t) =  m(x(t)) $ is also a solution of \eqref{eq:xisy}.
Therefore, the contraction property implies that
all the trajectories of \eqref{eq:nonlinsys} converge to the manifold $ \mathcal{M}:= \{x \in \R^n ~\vert~ p(x) = m(x) \}$.
\end{Remark}
%%  based on a similar argument to \cite[Lemma. 2]{sontag_cotraction_tutorial}

Let $\delta \xi_0 \in T \Omega_\xi$ denote an infinitesimal virtual variation to the initial condition $\xi(0) = \xi_0 \in \R^{\ell}$ of \eqref{eq:xisy}. Define
$
\delta \xi(t) := \frac{\partial \xi(t, \xi_0)}{\partial \xi_0}\delta \xi_0.
$
This yields the following variational system associated with \eqref{eq:xisy}:
\be \label{eq:varxi}
 \delta \dot{\xi}(t) =J_\xi(t, \xi, x) \delta \xi(t),
\ee
where $J_\xi(t, \xi, x) : = \frac{\partial f_\xi  }{\partial \xi}(t, \xi, x)$. This variational system will be instrumental in the subsequent analysis.

We can now state  our main results.

 %%%%%%%%%%%%%%%%%%%%

%%%%%%%%%%%%%%%%%%%%%%%%%%%%%%%%%%%%%%%%%%%%

%%%%%%%%%%%%%%%%%%%%%%%%%%%%%%%%%%%%%%%

%\section{Main results} 
%%%%%%%%%%%%%%%%%%%%%%%%%%%%%%%%%%%%%%%%%%%

\section{From partial contraction to horizontal contraction} \label{sec:main1}
%%%%%%%%%%%%%%%%%%%

\begin{Theorem} \label{prop:partohori}
Suppose that there exists $\lambda > 0$, and a Finsler-Lyapunov function $V_\xi : T \Omega_\xi \to \R_+$ for the system~\eqref{eq:xisy} such that 
	\be \label{eq:dvxi}
	\frac{\partial V_{\xi}(\xi, \delta \xi)}{\partial \xi} f_\xi (t, \xi, x) + \frac{\partial V_{\xi}(\xi, \delta \xi)}{\partial \delta \xi}J_\xi(t, \xi, x) \delta \xi \leq - \lambda  V_{\xi}(\xi, \delta \xi),
	\ee
for all $t \in \R_+$, $x\in \Omega$, and $(\xi, \delta \xi) \in T \Omega_\xi$. That is, the system~\eqref{eq:nonlinsys} is partially contractive w.r.t. $p(x)$. 
	%%%%%%%
	Consider a distribution $\mathcal{H}_x$ as in \eqref{eq:hqx} and the associated matrix $H(x) \in \R^{n \times \ell}$ in \eqref{eq:hqmatrix}. 
	%%%%%%%%%%
	Assume that there exists a vector norm $\vert \cdot \vert$, and~$d_6, d_7 >0$ such that
	\be \label{eq:normh}
	 d_6 \vert y \vert \leq \vert H^T(x) H(x) y \vert \leq d_7 \vert y \vert, \text{ for all } x\in \Omega,~y \in \R^\ell.
	\ee
    Furthermore, for all $t \in \R_+$ and all $x \in \Omega$,  assume that
	\be \label{eq:hfjxi}
	H_f^T(x) +H^T(x)J(t,x) = J_\xi (t, p(x), x) H^T(x),
	\ee
	where 
	\[
	H_f(x) : = \begin{bmatrix} \frac{\partial h_1(x)}{\partial x}f(t,x) & \cdots & \frac{\partial h_\ell(x)}{\partial x}f(t,x) \end{bmatrix} \in \R^{n \times \ell} 
	\]
	is the 
  the directional derivative of $H(x)$ along $f(t,x)$.
	Then, the system  \eqref{eq:nonlinsys} is also horizontally contractive w.r.t. $\mathcal{H}_x$.
\end{Theorem}
 
\begin{Remark} \label{re:pxhxmx}
It is usually difficult to find the matrix $H(x)$ in Theorem~\ref{prop:partohori} except some special cases. For example, consider the case in Remark~\ref{re:p-m}. Note that the manifold $\mathcal{M}$ is attractive, therefore we may conjecture that \eqref{eq:nonlinsys} is horizontally contractive along the directions which are orthogonal to $\mathcal{M}$. Hence, a possible choice of $H(x)$ is $H(x) = \frac{\partial^T p(x)}{\partial x} - \frac{\partial^T m(x)}{\partial x}$.
\end{Remark}

\begin{IEEEproof} [Proof of Thm.~\ref{prop:partohori}]
 Differentiating~$H^T(x) \delta x(t)$ along the trajectories of \eqref{eq:nonlinsys} and \eqref{eq:vardx} gives
	\be \label{eq:dhdx} \begin{aligned}
	\frac{\diff}{\diff t}(H^T(x)\delta x) =& (H_f^T (x) + H^T(x)J(t,x)) \delta x \\
	=& J_\xi (t, p(x), x) H^T(x) \delta x,
	\end{aligned} \ee
where \eqref{eq:hfjxi} is used. 
This implies that $\delta \xi(t) := H^T(x(t))\delta x(t)$ is a trajectory of \eqref{eq:varxi}. 
Recall that $\xi(t) = p(x(t))$ is a solution of \eqref{eq:xisy}.
	%%%%%%%%%%%%%%%%%%%
Let
	\be \label{eq:vbar}
	V(x, \delta x) := V_{\xi}(p(x), H^T(x)\delta x(t)).
	\ee
By \eqref{eq:hty},  $V(x, \delta x) = V(x, H(x) \delta x_h)$.
   Eq. \eqref{eq:vbar} implies that the derivative of $V(x, \delta x)$ is equal to the derivative of $V_{\xi}(\xi, \delta \xi)$ along the trajectories $\xi(t) = p(x)$ and $\delta \xi(t) = H^T(x)\delta x(t)$, i.e.,
   \[
    \frac{\diff}{\diff t} V(x(t), \delta x(t))= 
     \frac{\diff}{\diff t} V_\xi(\xi(t), \delta \xi(t))\vert_{ \xi = p(x),~\delta \xi = H^T(x)\delta x}.
   \]
   Specifically,
   \be \label{eq:dvbar} \begin{aligned}
      &\frac{\partial V (x, \delta x)}{\partial x} f(t,  x) + \frac{\partial V(x, \delta x)}{\partial \delta x}J(t,  x) \\
   =&  \bigg( \frac{\partial V_{\xi}(\xi, \delta \xi)}{\partial \xi} f_\xi (t, \xi, x) \\
   & + \frac{\partial V_{\xi}(\xi, \delta \xi)}{\partial \delta \xi}J_\xi(t, \xi, x) \delta \xi  \bigg)\vert_{ \xi(t) = p(x),~\delta \xi(t) = H^T(x)\delta x(t)} \\
   \leq& - \lambda V_\xi(\xi, \delta \xi)\vert_{ \xi(t) = p(x),~\delta \xi(t) = H^T(x)\delta x(t)}\\
   =& -\lambda  V(x, \delta x).
   \end{aligned}\ee	
In order to further show that $V(x(t), \delta x(t))$ is indeed a horizontal Finsler-Lyapunov function for~\eqref{eq:nonlinsys} w.r.t. $\mathcal{H}_x$, we need to find an associated function $F(x, \delta x)$ satisfying Definition~\ref{def:semifins}.

 Since $V_\xi(\xi, \delta \xi)$ is Finsler-Lyapunov function of the system~\eqref{eq:xisy}, there exists an associated function $F_\xi: T \Omega_\xi \to \R_+$ that satisfies Definition~\ref{def:semifins} with a zero vertical distribution. Let 
 \[
 F(x, \delta x) =: F_\xi(p(x), H^T(x) \delta x).
 \]
 Hence, 
 \[
d_4 \vert H^T(x) \delta x \vert \leq F(x,\delta x) \leq d_5 \vert H^T(x) \delta x \vert,
\]
for some $d_4, d_5 >0$. 
%%%%%%%%%%
By \eqref{eq:hty}, this yields
\[
d_4 \vert H^T(x) H(x) \delta x_h \vert \leq F(x,\delta x) \leq d_5 \vert H^T(x)  H(x) \delta x \vert.
\]
From \eqref{eq:normh}, we have
\[
d_4 d_6  \vert  \delta x_h \vert \leq F(x,\delta x) \leq d_5 d_7 \vert  \delta x_h \vert.
\]
That is, $F(x, \delta x)$ satisfies property \eqref{pp:d45} in Definition~\ref{def:semifins}. 
Furthermore, note that $V(x, \delta x)$ and $F(x, \delta x)$ satisfy condition~\eqref{eq:d123} and properties \eqref{pp:h2}-\eqref{pp:h5} in Definition~\ref{def:semifins} by directly inheriting those from $V_\xi(\xi, \delta \xi)$ and  $F_\xi(\xi, \delta \xi)$.
%%%%
Hence, $V(x, \delta x)$ is indeed a horizontal Finsler-Lyapunov function for \eqref{eq:nonlinsys}. By Prop.~\ref{prop:hori}, Eq.~\eqref{eq:dvbar} ensures that the system \eqref{eq:nonlinsys} is horizontally contractive w.r.t. $\mathcal{H}_x$. 
\end{IEEEproof}

The next example shows that condition \eqref{eq:hfjxi} generally does not hold even when $p(x)$ is linear.

\begin{Example}[Convergence to flow-invariant subspaces \cite{Pham2007}] \label{exa:1}
Let $H \in \R^{n \times \ell}$ and $Q \in \R^{n \times (n -\ell)}$, such that
\be  \label{eq:hqi}\begin{aligned}
& H^T H = I_\ell, \quad  Q^T Q  = I_{n -\ell}, \\
& H^T Q = 0,\quad H H^T + Q Q^T = I_n.
\end{aligned} \ee
Note that the above conditions hold if the column vectors of $H$ and $Q$, i.e., $h^1, \dots, h^\ell$, $q^1, \dots, q^{n-\ell}$, are \emph{orthonormal}, which means that they all have Euclidean norm one and are mutually orthogonal. Let $\mathcal{H}$ and $\mathcal{Q}$ denote the column subspaces of $H$ and $Q$, respectively. It is assumed that 
\be \label{eq:fuu}
f(t, \mathcal{Q} ) \subseteq \mathcal{Q},
\ee 
That is, $\mathcal{Q}$ is flow-invariant.
%%%%%%%%%%%%%%%%%%%%%%%%%%%%%%%%%
In this case, Eq. \eqref{eq:pxf} naturally holds with $p(x) = H^T x$ and
\be
g(t, p(x), x) := f(t, H p(x) + Q Q^T x  ).
\ee
Then, the system \eqref{eq:xisy} can be defined as
\be \label{eq:xisy1}
\dot \xi = f_\xi(t, \xi, x):= H^T g(t, H \xi + Q Q^T x ).
\ee
Note that $\xi(t) = H^T x(t)$ and $\xi(t) = 0$ are two particular solutions to \eqref{eq:xisy1}. Therefore, partial contraction w.r.t. $H^T x$ ensures that all the trajectories of \eqref{eq:nonlinsys} converge to $\mathcal{Q}$.
%%%%%%%%%%%%%%%

Then, we study if the above conditions also imply horizontal contraction. Consider the case when $f(t, x) = 0$ for all $x \in \mathcal{Q}$, which ensures that \eqref{eq:fuu} holds. In this case, the NTV system \eqref{eq:nonlinsys} can only be horizontally contractive w.r.t. $\mathcal{H}$. Without loss of generality, we can choose $H(x)$ in Theorem~\ref{prop:partohori} as $H$ here.
Note that $J_\xi(t, H^T x, x) = H^T J(t,x) H$. 
Then, condition \eqref{eq:hfjxi} in this case boils down to 
\be \label{eq:hjhh}
H^T J(t,x) = H^T J(t,x) H H^T,
\ee
which does not hold in general.
According to Theorem~\ref{prop:partohori}, if the system \eqref{eq:nonlinsys} is partially contractive w.r.t. $H^T x$, then \eqref{eq:hjhh} ensures that it is also horizontally contractive w.r.t. $\mathcal{H}$.
%%%%%%%%%%%%%%%%
From \eqref{eq:hqi}, a sufficient condition to ensure \eqref{eq:hjhh} is 
\be \label{eq:hjq}
H^T J(t,x) Q = 0.
\ee 
Consider a linear system, i.e., $J(t,x) = A \in \R^{n \times n}$. Then, Eq. \eqref{eq:hjq} holds naturally since \eqref{eq:fuu} implies $A \mathcal{Q} \subseteq \mathcal{Q}$.
\end{Example}
 
In general, for a horizontally contractive system \eqref{eq:nonlinsys}, the existence of $p(x)$ as in Definition~\ref{def:parcon}, which depends on the integrability of $H(x)$ (see e.g., Remark~\ref{re:pxhxmx}), is not ensured. That is, horizontal contraction does not necessarily imply partial contraction. 

\section{From horizontal contraction to $k$-contraction}  \label{sec:main2}

To facilitate the subsequent result, we first prove a useful property of wedge products. Recall that a vector norm~$|\cdot|:\R^n\to\R_+$ is called monotonic if for any~$x,y\in\R^n$ with~$|x_i|\leq |y_i|$ for all~$i$, we have~$|x|\leq|y|$~\cite{mono_norms}. For example, all the $L_p$ norms are monotonic. 
%%%%

\begin{Lemma} \label{le:wedge}
	Pick $k \in [2,n]$. Consider a set of $k$ time-varying vectors $a^1(t), \cdots, a^{k}(t) \in \R^n$. Assume that there exist constants $\ell \in [1, k-1]$, $\gamma_1 \geq 1 $, $\gamma_2 \geq 1, \beta >0$,  and a monotonic vector norm $\vert \cdot \vert$ such that
	\begin{align}\label{eq:kpara} 
	 \vert a^j(t) \vert& \leq \gamma_1 \exp(-\beta t)  \vert a^j(0) \vert,  &j&= 1, \dots, \ell, \\
	 \vert a^{\ell +i}(t) \vert& \leq \gamma_2 \vert a^{\ell +i}(0) \vert , &i& =1, \dots, k -\ell,
	\end{align}
	 for all $ t \in \R_+$. 
Then, $|\wedge_{j=1}^k a^j(t)|$ decays to zero exponentially. Furthermore, for $|\wedge_{j=1}^k a^j(0) | \neq 0$, there exists $\bar{\gamma}, \bar{\beta} >0$ such that 
\be
  |\wedge_{j=1}^k a^j(t)| \leq \bar{\gamma} \exp(- \bar{\beta} t) |\wedge_{j=1}^k a^j(0) |, \text{ for all } t \in \R_+.
\ee
\end{Lemma}

\begin{IEEEproof}
	Recall that $ \vert \wedge_{j=1}^k a^j(t) \vert$ is the $k$-volume of the $k$-parallelotope generated by $a^1(t), \cdots, a^k(t)$. Intuitively speaking, \eqref{eq:kpara} implies that at least one edge of this $k$-parallelotope shrinks exponentially, and the other edges are uniformly bounded. Therefore, its $k$-volume also shrinks exponentially. Here, we only provide a detailed proof for $\ell = k-1$. The proof for other cases is based on similar arguments.
	
	By the property of the wedge product, we have
	\be \label{eq:recur}\begin{aligned}
	%%%
	\wedge_{j=1}^k a^j 
	=(\wedge_{j=1}^{k-1}a^j)  \wedge \sum_{i=1}^n a^k_i e^i
	=\sum_{i=1}^n a^k_i     (\wedge_{j=1}^{k-1}a^j) \wedge e^i.
	%%%%
	\end{aligned}\ee
	where $a_i^k$ denotes the $i$th entry of $a^k$, and $e^i$ denotes the $i$th canonical vector in $\R^n$.
	Let
	\[z^i(t):= (\wedge_{j=1}^{k-1}a^j(t)) \wedge e^i=\begin{bmatrix}  a^1 &\dots& a^{k-1}& e^i \end{bmatrix}^{(k)}.\]
	By the Leibniz formula for determinants, $z^i(t)$ has at least~$\binom{n-1}{k}$ 
	zero entries, and every nonzero entry is an entry of the vector~$\wedge_{j=1}^{k-1}a^j$
	multiplied by either plus one or minus one. Hence, for any monotonic vector norm~$|\cdot|:\R^n\to\R_+$, we have~$|z^i(t)|\leq |\wedge_{j=1}^{k-1}a^j(t)|$ and thus
	\[\label{eq:pou1} \begin{aligned}
	|\wedge_{j=1}^k a^j(t) | \leq &
	\sum_{i=1}^n  (|a^k_i|      |\wedge_{j=1}^{k-1}a^j(t)|) \\
	\leq &  n \gamma_2  |a^k(0) |
	 |\wedge_{j=1}^{k-1}a^j(t)|,
	\end{aligned}\]
	where we used the fact $|a^k_i(t)| \leq |a^k(t)| \leq \gamma_2 |a^k(0)| $. Then, by repetitively using the recursive formula~\eqref{eq:recur}, 
	we have
	\[
	|\wedge_{j=1}^k a^j(t)| \leq (n\gamma_1) ^{k-1} \gamma_2  \exp(- (k-1)\beta t)  \prod_{j=1}^k \vert a^j(0) \vert.
    \]
    That is, $|\wedge_{j=1}^k a^j(t)|$ converges to zero exponentially.
     For the case $|\wedge_{j=1}^k a^j(0) | \neq 0$, note that there exists a large enough constant $\gamma_3 >0$ such that $\prod_{j=1}^k \vert a^j(0) \vert \leq   \gamma_3	|\wedge_{j=1}^k a^j(0) |  $. Hence,  
     	\[
     |\wedge_{j=1}^k a^j(t)| \leq \bar{\gamma} \exp(- \bar{\beta} t) |\wedge_{j=1}^k a^j(0) |,
     \]
     where $\bar{\gamma} :=  (n \gamma_1)^{k-1} \gamma_2 \gamma_3  $ and $\bar{\beta} := (k-1)\beta$.
\end{IEEEproof}

The next result specifies a sufficient condition such that horizontal contraction implies $k$-contraction. In the subsequent analysis, we assume that the vector norms under concern are monotonic. 

\begin{Theorem} \label{prop:horitok}
	Suppose that the system \eqref{eq:nonlinsys} 
	is horizontally contractive w.r.t. $\mathcal{H}_x$ defined in \eqref{eq:hqx}, where we can write $\ell = n-k+1$ for some $k = [1, n]$. 
	%%%%%%%%%%%%%%%%%%%
	For the matrices $H(x) \in \R^{n \times (n-k+1)}$ and $Q(x) \in \R^{n \times (k-1)}$ given in \eqref{eq:hty}, define 
	\be \label{eq:defm}
	M(x) := \begin{bmatrix} 
		H^T(x) \\
		Q^T(x)
	\end{bmatrix} \in \R^{n \times n}.
	\ee
    Assume that for some vector norm $\vert \cdot \vert$, there exist $c_i>0$, $i=1,\dots,6$, such that for all $x \in \Omega$
    \be \label{eq:normhq} \begin{aligned} 		 
   & c_1 \vert y \vert \leq \vert H^T(x) H(x) y \vert \leq c_2 \vert y \vert, ~y \in \R^{n-k+1}, \\
     & c_3 \vert y \vert \leq \vert M(x) y \vert \leq c_4 \vert y \vert,~y \in \R^n, \\
     & c_5 \vert y \vert \leq \vert M^{(k)}(x) y \vert \leq c_6 \vert y \vert, ~y \in \R^{\binom{n}{k}}.
    \end{aligned}\ee
	%%%%%%%%%%%%%%
    Furthermore, for any two trajectories of the system \eqref{eq:nonlinsys}, denoted $x^1(t)$ and $x^2(t)$, there exists a constant $\gamma_1 > 1$ such that 
	\be \label{eq:unfbd}
	|x^1(t) - x^2(t)| \leq  \gamma_1 |x^1(0) - x^2(0)|, \quad \text{for all } t \geq 0.
	\ee
  Then, the system \eqref{eq:nonlinsys} is also $k$-contractive on $\Omega$. 
\end{Theorem}

\begin{IEEEproof}
%%%%%%%%%%%%%
First, we show that the variational system \eqref{eq:vardx} is uniformly bounded under the condition~\eqref{eq:unfbd}. 	
%%%%%%%%%%%%%%%%%%
Pick $\varepsilon >0$, and $x_0, z_0 \in \R^n$.  Define
\begin{equation}
z(t, \varepsilon) := \frac{x(t, x_0 + \varepsilon z_0 ) - x(t, x_0)}{\varepsilon}.
\end{equation}
By \eqref{eq:unfbd}, we have
\begin{equation} \label{eq:xtxi0}
\vert z(t, \varepsilon) \vert \leq \gamma_1 | z_0|  .
\end{equation}
Since the solutions of \eqref{eq:nonlinsys} are continuously dependent on the initial conditions, the limit $z(t) := \lim_{\varepsilon \to 0} z(t, \varepsilon)$ exists and 
$
z(t) = \frac{\partial x(t, x_0)}{\partial x_0} z_0.
$
This implies that $z(t)$ is a solution of the variational system \eqref{eq:vardx} with the initial condition $z(0) = z_0$.
Thus, \eqref{eq:xtxi0} reduces to
\be \label{eq:xiuni}
\vert z(t) \vert \leq \gamma_1 | z_0|.
\ee
That is, the variational system \eqref{eq:vardx} is uniformly bounded. 
%%%%%%%%%%%%

Since the system \eqref{eq:nonlinsys} is horizontally contractive, Definition~\ref{def:hori} implies that there exist  $\gamma_2 > 1$ and $\beta> 0$ such that
\be \label{eq:ceta}
\vert \delta x_h(t) \vert \leq  \gamma_2 e^{-\beta t}  \vert \delta x_h(0) \vert, \quad \text{for all } t\geq 0.
\ee

Recall that any $\delta x \in T_x \Omega$ can be rewritten as $\delta x  =  H(x) \delta x_h + Q(x) \delta_q$, where $\delta x_h \in \R^{n-k+1}$, and $\delta x_q \in \R^{k-1}$. Hence, 
\be
M(x)\delta x  = 
 \begin{bmatrix} 
	H^T(x) H(x) \delta x_h \\
	Q^T(x) Q(x) \delta x_q
 \end{bmatrix},
\ee
where we use the fact that $H^T(x) Q(x) = 0$.
From \eqref{eq:normhq} and \eqref{eq:ceta}, we have
\be \label{eq:n-k+1} \begin{aligned}
\vert H^T(x(t)) H(x(t)) \delta x_h(t) \vert \leq &  c_2 \vert \delta x_h(t) \vert  \\
 \leq&  c_2 \gamma_2 e^{-\beta_1 t}  \vert \delta x_h(0) \vert,
\end{aligned}\ee
for all $t \geq 0$.
That is, the first $(n-k+1)$ entries of $M(x)\delta x$ converge to zero exponentially. 
%%%%%%%%%%%%
By virtue of monotonic vector norms, 
\be \label{eq:k-1} \begin{aligned}
\vert 	Q^T(x(t)) Q(x(t)) \delta x_q(t) \vert \leq&
 \vert M(x(t)) \delta x(t) \vert  \\
 \leq & c_4 \vert \delta x(t) \vert \\
 \leq & c_4 \gamma_1 | \delta x(0)|  , \text{ for all } t \geq 0,
\end{aligned}\ee
where \eqref{eq:normhq} and \eqref{eq:xiuni} are used. 
That is, the last $(k-1)$ entries of $M(x)\delta x$ are uniformly bounded.

Then, as defined in \eqref{eq:defy}, consider $k$ trajectories of the variation system \eqref{eq:vardx} specified by $x(t, a)$ with $a \in \Omega$, i.e., $\delta^1 x(t, a), \dots, \delta^k x(t, a),$ such that $\wedge_{i=1}^k \delta^i a \neq 0 $. Let 
\[
A(t) := 
\begin{bmatrix}
	M(x(t, a)) \delta x^1(t, a) & \cdots &  M(x(t, a)) \delta x^k(t, a)
\end{bmatrix}^T.
\]
Note that $A \in \R^{k \times n}$. 
Using \eqref{eq:prodab} and \eqref{eq:wedgdef}, we have
\begin{equation}\begin{aligned}
 (A^T(t))^{(k)} = & \wedge_{i=1}^k M(x(t,a)) \delta x^i(t, a)\\
       =&  M^{(k)}(x(t, a)) y(t,a) ,
\end{aligned}\end{equation}
where $y(t, a ) := \wedge_{i=1}^k \delta x^i(t,a)$ as in \eqref{eq:defy}.
%%%%%%%%%%%%%%%%
Let $a^i(t)$ denote the $i$th column of $A(t)$. From \eqref{eq:n-k+1} and \eqref{eq:k-1}, $\vert a^i(t) \vert$ with $i =1, \dots, n-k+1$, converges to zero exponentially, and $\vert a^i(t) \vert$ with $i =n-k+2, \dots, n$, are uniformly bounded. 

Let~$Q_{k,n}$ denote the set of increasing sequences of~$k$ numbers from~$[1,n]$
ordered lexicographically. With a slight abuse of notation, we treat such ordered sequences as sets. The cardinality of $Q_{k,n}$ is $r:=\binom{n}{k}$. The $j$th element of $Q_{k,n}$ is denoted $\kappa_j$. For example,
$
Q_{2,3} =\{\kappa_1, \kappa_2, \kappa_3\}, \text{ with } ~\kappa_1= \{1,2\},~\kappa_2 = \{1,3\}, ~\kappa_3 = \{2,3\}.
$
Then, $\wedge_{i \in  \kappa_j} a^i(t)$ is the $j$th entry of $M^{(k)}(x(t, a)) y(t,a)$. Note that any $k$ vectors $a^i(t)$, $i \in \kappa_j$, satisfy condition \eqref{eq:kpara} in Lemma~\ref{le:wedge}. Hence, it is ensured that $\wedge_{i \in  \kappa_j} a^i(t)$, $j = 1, \cdots, r$, i.e., all the entries of $M^{(k)}(x(t, a)) y(t,a)$, converge to zero exponentially. Therefore, there exist $\gamma_3>1 $ and $ \eta >0$ such that 
\be
\vert M^{(k)}(x(t, a)) y(t,a) \vert \leq  \gamma_3 e^{-\eta t} \vert M^{(k)}(x(0, a)) y(0,a) \vert,
\ee
for all $t \geq 0$. Then from \eqref{eq:normhq}, we have 
\be
\vert  y(t,a) \vert \leq  \frac{c_6 \gamma_3 }{c_5} \exp(-\eta t) \vert  y(0,a)  ,
\ee
 for all $t \geq 0$. This implies that the system \eqref{eq:nonlinsys} is $k$-contractive by Definition~\ref{def:kord}.
\end{IEEEproof}

\begin{Remark} \label{re:unfbd}
Condition~\eqref{eq:unfbd} actually asserts uniform incremental stability, that is, every solution of \eqref{eq:nonlinsys} is uniformly stable. By Coppel’s inequality (see e.g., \cite{vidyasagar2002nonlinear}), condition~\eqref{eq:unfbd} is ensured if there exists a matrix measure $\mu: \R^{+} \times \R^n \to \R$ and a constant $c>0$ such that 
\[
\int_0^t \mu(J(s,x(s,a))) \diff s \leq c, \text{ for all } a \in \Omega, t \geq 0.
\]
As shown in the proof of Theorem~\ref{prop:horitok}, condition~\eqref{eq:unfbd} can actually be relaxed since it is only required that the exponential convergence rate of $H^T(x) H(x)\delta x_h $ is larger than the expansion rate of $Q^T(x) Q(x)\delta x_q $ to conclude that $ M(x) \delta x $ converges to zero exponentially. For example, consider the  linear system
\begin{align*}
    \dot x_1 &= c_1 x_1 + x_2, \\
    \dot x_2 &= - c_2 x_2,
\end{align*}
where $c_2 > c_1 >0 $. This system is both $2$-contractive and horizontally contractive w.r.t. $\mathcal{H}_x := \spanop\{ \begin{bmatrix} 0 & 1 \end{bmatrix}^T \}$. However, it is not uniformly stable since $x_1(t)$ will go to infinity. Note also that this system is still horizontally contractive but not $2$-contractive with $c_1 > c_2 >0 $.   
\end{Remark}

%%%%%%%%%%%%%%%%
Since $k$-contraction does not require the existence of the distributions $\mathcal{H}_x$ and $\mathcal{Q}_x$. Therefore, $k$-contraction does not ensure horizontal contraction in general. This is shown by the next example.

\begin{figure}[t]
	\begin{center}
		\centering
		\includegraphics[scale=0.46]{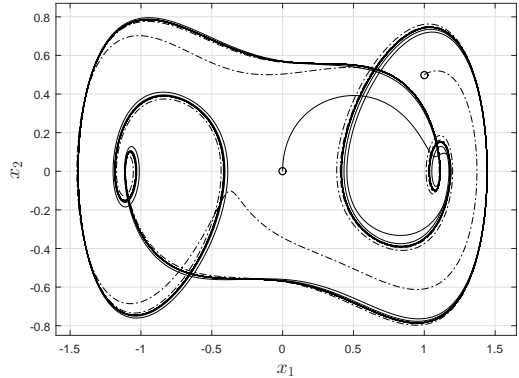} 
		\caption{Trajectories of the Duffing oscillator for two different initial conditions (black circle). }
		\label{duffing} 
	\end{center}
\end{figure}

\begin{Example}[The forced Duffing oscillator \cite{thompson2002nonlinear}]
	Consider the NTV system 
	\begin{equation} \label{eq:duffing}\begin{aligned}
	& \dot{x}_1 = x_2, \\
	& \dot{x}_2 = \theta_1 x_1 - \theta_2 x_1^3 - \theta_3 x_2 + \theta_4 \cos(\theta_5 t)  ,
	\end{aligned}\end{equation}
	where $\theta_i > 0$, $i = 1, \dots, 5$. For $\theta_1 = \theta_2 = 1, \theta_3 = 0.3, \theta_4 = 0.37$ and $\theta_5 = 1.2$, this system has an attractor with self intersections as shown in Fig. \ref{duffing}. Therefore, this system is not contracting in any specified direction, that is, it is not horizontally contractive. Indeed, it can have more complicated attractors with other different parameters. However, the Jacobian~$J(t,x)$ of~\eqref{eq:duffing},   satisfies $J^{[2]}(t,x) = -\theta_3 <0$. Therefore, it is $2$-contractive according to Definition~\ref{def:kord} and \cite[Thm. 4]{kordercont}.
\end{Example}

\section{An Example: the Andronov-Hopf oscillator} \label{sec:aho}
In this section, the Andronov-Hopf oscillator is revisited to illustrate our results. It is shown that this system is partially contractive, horizontally contractive, and $2$-contractive. 

The dynamics of	the Andronov-Hopf oscillator  is	%%%%%%%%%%%%%%%%%%%%%%%%%%%%%%%%%%%%%%%%%%%%%%%%%%%%%
		\begin{align}\label{eq:unitcir} 
		\dot{x}_1 &= -x_2 - x_1(x_1^2 + x_2^2 -1), \nonumber \\
		\dot{x}_2& = x_1 - x_2(x_1^2 + x_2^2 - 1).
		\end{align} 
		%%%%%%
Note that~$x = 0$ is an unstable equilibrium, and the unit circle is a stable limit cycle.
%%%%%%%%
The associated variational system is 
	\begin{equation} \label{eq:delta_x}
		\delta \dot x(t) = J(x(t)) \delta x(t), 
	\end{equation}
	where
	\begin{equation}
		J(x)  = 
		\begin{bmatrix}
		1-3x_1^2-x_2^2 & -2x_1 x_2 -1 \\
		-2x_1 x_2+1  & 1-x_1^2-3x_2^2  
		\end{bmatrix}.
\end{equation}
		%%%%%%
		
\begin{Proposition} 
	Fix $0 < \gamma_1 < 1 < \gamma_2$. Consider the manifold $\Omega: =\{ x\in \R^2 ~\vert ~ \gamma_1 \leq x_1^2 + x_2^2 \leq  \gamma_2 \}$. The system~\eqref{eq:unitcir} is partially contractive, horizontally contractive, and $2$-contractive on $\Omega$.
\end{Proposition}

\begin{IEEEproof} Note that $\Omega$ here is forward invariant and connected, but not convex.
%%%%%%%%%%%%%%%%%%	
	
	(i) \emph{Partial contraction}: 
	 Let $p(x) :=  x_1^2 + x_2^2 -1$, and 
	\be
	g(\xi, x) :=
	\begin{bmatrix}
		-x_2 -x_1 \xi \\
		x_1 -x_2 \xi
	\end{bmatrix}.
	\ee
	Then $g(p(x) , x) = f(x)$. That is, condition \eqref{eq:pxf} in Definition~\ref{def:parcon} is satisfied. In this case, the system \eqref{eq:xisy} reduces to 
	\be \label{eq:xisys}
	\dot{\xi}(t) = -2(x_1^2(t) + x_2^2(t)) \xi(t).
	\ee 
	Clearly, \eqref{eq:xisys} is contractive for all $ x \in \Omega$.
	In this case, a Finsler-Lypapunov function for \eqref{eq:xisys} can be selected as $V_\xi(\xi, \delta \xi) = \delta \xi^2$, and we have 
	$\frac{\partial V_\xi(\xi, \delta \xi) }{ \partial \delta \xi} \leq -4 \gamma_1 V_\xi(\xi, \delta \xi) $ since $ \delta \dot{\xi}(t) = -2(x_1^2(t) + x_2^2(t)) \delta \xi(t)$.
	 Therefore, the system \eqref{eq:unitcir} is partially contractive w.r.t. $p(x) = x_1^2 + x_2^2 -1$. Note that both $p(x)$ and the origin are the solutions to \eqref{eq:xisys}. Therefore, for any $a \in \Omega$, 
	\be \label{eq:pexpon}
	\vert p(x(t,a))  \vert \leq e^{-2 \gamma_1 t} \vert p(a) \vert, \text{ for all } t \in \R_+.
	\ee
	%%%%%%%%%%%%%%%%%%%
	
	(ii) \emph{From partial contraction to horizontal contraction}: Let 
	\be
	H(x) = \frac{1}{x_1^2 + x_2^2} \begin{bmatrix}  x_1 \\  x_2 \end{bmatrix}.
	\ee
	Note that condition \eqref{eq:normh} in Theorem~\ref{prop:partohori} holds since
	$
	  \gamma_2^{-1}  \leq  H^T(x) H(x) = (x_1^2 + x_2^2)^{-1} < \gamma_1^{-1}, \text{ for all } x\in \Omega
	$. 
	Based on a straightforward calculation, 
	\[
	\frac{\diff}{\diff t} (H^T(x(t)) \delta x(t)) = -2 (x_1^2(t) + x_2^2(t))H^T(x(t)) \delta x(t) ,
	\]
	which implies that condition \eqref{eq:hfjxi} in Theorem~\ref{prop:partohori} also holds. Therefore, Theorem~\ref{prop:partohori} implies that this system is also horizontally contractive w.r.t. 
	$
	\mathcal{H}_x :=\spanop \left \{\frac{1}{x_1^2 + x_2^2} \begin{bmatrix}  x_1 \\  x_2 \end{bmatrix} \right \}. 
	$
	Furthermore, as shown in the proof of Theorem~\ref{prop:partohori}, a horizontal Finsler-Lyapunov function can be constructed as:
	\[\begin{aligned}
	V(x, \delta x) =& V_\xi(p(x), H^T(x)\delta x)  \\
	     =&  \delta^T x H(x) H^T(x)\delta x \\
	     =&  \left (\frac{x_1 \delta x_1 + x_2 \delta x_2}{x_1^2 + x_2^2}  \right)^2.
	\end{aligned}\]
	Indeed, $\dot{V}(x, \delta x) = -4(x_1^2 + x_2^2)V(x, \delta x) \leq -4 \gamma_1 V(x, \delta x)$ for all $x \in \Omega$.
	%%%%%%%%%%%%%%%%%%%%%%%%%%%%%%%%
	
  (iii) \emph{From horizontal contraction to $2$-contraction}: The matrix~$
  J_{sym} (x) := (J(x) + J^T(x))/{2} 
 $ has eigenvalues: $\lambda_1 (x)= 1 - x_1^2 -x_2^2$,  $\lambda_2(x) = 1 - 3 x_1^2 - 3x_2^2$. Hence, 
  \[
  \mu_2(J (x)) = 1 - x_1^2 -x_2^2,
  \]
  where $\mu_2(\cdot)$ is the matrix measure associated with the $L_2$ norm. By~\eqref{eq:pexpon}, we have 
  \[
  \int_0^t \mu_2(J (x(s,a))) \diff s \leq \frac{1}{2 \gamma_1}\max\{\gamma_2 -1, 1-\gamma_1\},
  \]
  for all $a \in \Omega$ and $t \geq 0$. This implies that condition~\eqref{eq:unfbd} holds as shown in Remark~\ref{re:unfbd}.
  %%%%%%%%%%%
    Let $H(x) := \frac{1}{x_1^2 + x_2^2} \begin{bmatrix} x_1 & x_2 \end{bmatrix}^T$, and $Q(x) := \frac{1}{x_1^2 + x_2^2} \begin{bmatrix} -x_2 & x_1 \end{bmatrix}^T$.  Condition \eqref{eq:normhq} holds for all $x \in \Omega$ since $\Omega $ is compact. Therefore, Theorem~\ref{prop:horitok} ensures that the system~\eqref{eq:unitcir} is $2$-contractive.
  %%%%%%%%%%
  
  In order to show this explicitly, for the $2$nd compound equation of \eqref{eq:delta_x}, i.e., $\dot y(t) = J^{[2]}(x(t)) y(t) $, we consider a change of coordinate $y_m = M^{(2)}(x) y = \frac{y}{x^2_1 + x_2^2}$ with  $M(x) := \begin{bmatrix} H(x)  & Q(x) \end{bmatrix}^T$.
  Then, 
  \begin{align*} \label{eq:ym} 
  \dot{y}_m(t) = &M^{(2)}(x) J^{[2]}(x) (M^{(2)}(x))^{-1}y_m(t) + \\      &M^{(2)}_f(x)  (M^{(2)}(x))^{-1} y_m(t) \\
     =& -2(x_1^2 + x_2^2) y_m(t),
  \end{align*}
  where $M^{(2)}_f(x)$ denotes the directional derivative of $M^{(2)}(x)$ along vector field of~\eqref{eq:unitcir}. This implies that $\vert y(t) \vert$ decays to zero exponentially, that is, the system~\eqref{eq:unitcir} is $2$-contractive on $\Omega$ according to Definition~\ref{def:kord}.
\end{IEEEproof}

 \section{Conclusions} \label{sec:con}
%%%%%%%%%%%%%%%%%%%%%%%%%%%%%%%%%%%%%%%%%%%%%%%%%%%%%%%
This note shows that partial contraction, horizontal contraction, and $k$-contraction are not equivalent in general. Some sufficient conditions are specified such that $k$-contraction is achieved from partial and horizontal contraction. Since it is known that partial and horizontal contraction can be used to solve synchronization problems of networked systems, a related research direction is to study how $k$-contraction can simplify the analysis of synchronization problems.

As shown by \cite{muldo1990} and a recent related work \cite{angeli2020robust}, $k$-contraction with $k=2$ is very effective to study nonlinear time-invariant systems, in particular, to rule out oscillatory behaviors, and then almost all trajectories converge to an equilibrium (not necessary unique). Therefore, 
 for some partially or horizontally contractive nonlinear time-invariant systems, if the derived conditions hold and lead to $2$-contraction, then we can draw strong conclusions regarding their global asymptotic stability. As a straightforward application, this may simplify the I$\&$I stabilization method in the sense that the target system there is not required to be specified. 
 
 \begin{comment}
 the condition (A1) in \cite[Prop. 1]{wang2016immersion}, i.e, the so-called target system specified in I$\&$I stabilization procedure has a globally asymptotically stable equilibrium, can be simply replaced by that the system in concern has a unique equilibrium. That is, there is no need to specify the target system in this case, which further simplifies the I$\&$I procedure.
%%%
As a particular application, the proposed results can be utilized to simplify the immersion $\&$ invariance stabilization procedure for certain specified scenarios.

Consider a nonlinear time-invariant system whose solutions evolve on a forward-invariant, closed, and convex manifold $\Omega$. Assume that it is $2$-contractive on $\Omega$. Then, as shown by \cite{muldo1990} and a recent related work \cite{angeli2020robust}, this system has no non-trivial periodic solutions, i.e., limit cycles.  Furthermore, if it has a unique equilibrium in $\Omega$. Then, every solution emanating from $\Omega$ converges to this equilibrium \cite{li1995}.
\end{comment}

%%%%%%%%%%%%%%%%%

\bibliographystyle{IEEEtranS}
\bibliography{part_contraction}

 \end{document}